\documentclass{SCIYA2017enOL}
\online
\usepackage[T1]{fontenc}
\usepackage{amssymb}

\makeatletter
\let\th@plain\th@newdefinition
\let\theorem\relax
\let\c@theorem\relax
\makeatother

\usepackage[nonotes]{degt}

\def\2{\color{red}}
\let\2\relax

\def\3{\color{red}}
\let\3\relax

\def\newcs#1#2{\expandafter\gdef\csname#1\endcsname{#2}}

\def\BB{$B2$}
\let\NMAX=N
\def\Nmax#1{\NMAX_{#1}}
\def\Nstar#1{\NMAX^*_{#1}}

\def\G{\frak G}
\def\bS{\bold S}
\def\lattice{\bS\sb{h}}
\def\Sperp{\bold{T}}
\def\(#1){[#1]}
\def\L{\bold L}
\def\th{\tilde h}

\def\Golay{\Cal{G}}
\def\CC{\Cal{C}}
\def\CK{\Cal{K}}
\def\CO{\Cal{O}}
\def\CS{\Cal{S}}
\def\KC{\bar\CO}
\def\KK{\bar\CK}
\def\KC{\CO_*}
\def\KK{\CK_*}
\let\dm\Omega
\let\graph\Gamma

\let\eq\sim
\let\Eq\approx
\def\cl#1{[#1]}
\def\Cl#1{[\![#1]\!]}
\def\sd{\mathbin\vartriangle}

\def\supp{\operatorname{supp}}
\def\Fn{\operatorname{Fn}}
\def\Km{\operatorname{Km}}
\def\Fano{\Cal{F}}
\def\irr{^*}

\def\EQ#1){}%
\let\NE\to

\def\OGplus{\OG^{+\!}}

\def\famdm#1{\tilde{#1}}

\def\famX{\Cal X}

\def\famXdm{\famdm\famX}

\def\famB{\Cal B}

\def\br{\famdm r}
\def\bc{\famdm c}

\def\gk{\frak k}
\def\go{\frak o}
\def\gr{\frak r}
\def\gs{\frak s}
\def\gu{\frak u}
\def\ggu{\gu'}
\def\gv{\frak v}
\def\ggv{\gv'}
\def\gw{\frak w}
\def\ggw{\gw'}

\def\hu{\hat u}

\let\hh\hbar
\def\ds#1#2{\mathopen\|#1/#2\mathclose\|}

\def\cn#1-#2{#1\sb{#2}}
\def\ln#1-#2{#1^\star\sb{#2}}

\def\dual{^\vee}

\def\CN#1{#1}
\def\LN#1{#1^\star}
\def\STRAT#1{\ifcase#1 \or
 \LN1\or
 \LN2\or
 \CN3\or
 \CN4\or
 \CN5\or
 \CN6\or
 \CN7\fi}
\def\astrat#1{\raise6pt\hbox{\hypertarget{strat-#1}{}}\STRAT#1}
\def\rstrat#1{\hyperlink{strat-#1}{\STRAT#1}}
\pdef\strat#1{\ifmmode\rstrat#1\else$\rstrat#1$\fi}

\def\rsame#1#2{\hyperlink{same-#1-#2}{${}^{#2}$}}
\def\same#1{\expandafter\ifx\csname same-\SAME-#1\endcsname\relax
 \smash{\llap{\raise9pt\hbox{\hypertarget{same-\SAME-#1}{}}${}^{#1}$\ }}
 \expandafter\gdef\csname same-\SAME-#1\endcsname{#1}
\else\smash{\llap{\rsame\SAME#1\ }}\fi}

\pdef\symplectic(#1,#2){\expandafter\ifx\csname(#1,#2)\endcsname\relax??\else
 \csname(#1,#2)\endcsname\fi}

\def\mrsign#1{$\ifcase#1\or*\or\dagger\fi$}
\def\maxr#1{\llap{\smash{\raise9pt\hbox{\hypertarget{maxreal-#1}{}}}\mrsign#1\ }}
\pdef\maxreal#1{\hyperlink{maxreal-#1}{\mrsign#1}}

\def\tabdefs{%
 \let\sep\
 \let\++%
 \def\(##1,##2){\symplectic(##1,##2)}%
 \def\^##1{\ifx1##1\else\rlap{$^##1$}\fi}%
 \let\s\strat
}

\begin{document}

\ensubject{fdsfd}


\ArticleType{ARTICLES}
\SpecialTopic{Partially supported by T\"{U}B\DOTaccent{I}TAK grant 118F413}
\Year{2023}
\Month{January}%
\Vol{66}
\No{1}
\BeginPage{1} %
\DOI{10.1007/s11425-016-5135-4}
\ReceiveDate{January 1, 2022}
\AcceptDate{January 1, 2022}
\OnlineDate{January 1, 2022}

\title[]{Conics on Barth--Bauer octics}{Conics on Barth--Bauer octics}

\author[1]{Alex Degtyarev}{degt@fen.bilkent.edu.tr}

\AuthorMark{Degtyarev A}

\AuthorCitation{Degtyarev A}

\address[1]{%
Department of Mathematics\\
Bilkent University\\
06800 Ankara, TURKEY}



\keywords{%
$K3$-surface, octic surface, Kummer surface, conic, Mukai group%
}

\MSC{%
Primary: 14J28;
Secondary: 14N25%
}


\abstract{%
We analyze the configurations of conics and lines on a special class of
Kummer octic surfaces. In particular, we bound the number of conics by $176$
and show that there is a unique surface with $176$ conics, all irreducible:
it admits a faithful action of one of the Mukai groups.
Therefore, we also discuss conics and lines
on Mukai surfaces: we discover a double plane
(ramified at a smooth sextic curve) that contains $8910$ smooth conics.
}

\maketitle

\section{Introduction}\label{S.intro}

The classical, going back to
A.~Cayley~\cite{Cayley:quartics,Cayley:cubics},
G.~Salmon~\cite{Salmon:1862},
F.~Schur~\cite{Schur:quartics}, and
B.~Segre~\cite{Segre}, problem of counting/estimating the number of smooth
rational curves on polarized algebraic surfaces has become increasingly
popular in the last decade or so. In spite of considerable efforts (\cf.
\cite{Barth.Bauer:conics,
Bauer:conics,
Boissiere.Sarti,
Caporaso,
rams.schuett:quintics}),
at present,
apart from the ``trivial'' cases of quadric and cubic surfaces,
only for \emph{lines} (smooth rational curves of projective degree~$1$) and only on
polarized \emph{$K3$-surfaces} a satisfactory answer is known, see
\cite{degt:lines,
degt:supersingular,
DIS,
rams.schuett:char3,
rams.schuett,
rams.schuett:char2,
Veniani} and further references therein.
In this paper, still working with polarized $K3$-surfaces, we make a step
towards understanding the maximal number of \emph{conics}, \ie, smooth
rational curves of degree~$2$.
Remarkably, in our approach we do not need to \emph{require} that the
conics should be smooth: it appears that the presence of singular (reducible)
ones reduces the upper bound, see \autoref{conj.N*}. This is yet another
mystery still to be understood.
For the moment, the sharp upper bound $\Nmax{2n}(2)$ on
the number of conics is known only for sextic $K3$-surfaces in~$\Cp4$: one
has $\Nmax6(2)=285$, see~\cite{degt:conics}.

Recall that the \emph{Kummer surface} $\Km(A)$ of an abelian surface~$A$ is
the quotient $A/{\pm1}$ blown up at the sixteen nodes---the images of the
sixteen fixed points of the involution. As is well known, $\Km(A)$ is a
$K3$-surface equipped with a distinguished collection of sixteen pairwise
disjoint smooth rational curves, \viz. the exceptional divisors contracted by
the projection $\Km(X)\to A/{\pm1}$. Conversely
(Nikulin~\cite{Nikulin:Kummer}), any $K3$-surface with sixteen pairwise
disjoint $(-2)$-curves is Kummer.

Extending the construction of Barth--Bauer~\cite{Barth.Bauer:conics} and
Bauer~\cite{Bauer:conics}, we define a \emph{Barth--Bauer surface} of
degree $h^2=2n\in2\Z^+$ as a
polarized Kummer surface
$X\into\Cp{n+1}$
{\2with the property} that the sixteen Kummer divisors map to sixteen
irreducible conics in $\Cp{n+1}$. Conjecturally
(see~\cite{degt:4Kummer,degt:800}), the maximal number of irreducible conics
on a smooth quartic surface is $\Nmax4(2)=800$, and this maximum is attained
at a certain Barth--Bauer quartic.
Therefore, in this paper we make an attempt to estimate the maximum
$\Nmax8(2)$ for octic surfaces
by obtaining a complete classification of the
Barth--Bauer octics up to \emph{equiconical deformation}, \ie, deformation in
$\Cp5$ preserving the bi-colored \emph{full Fano graph}
\[*
\Fn X:=\Fn_1X\cup\Fn\irr_2X
\]
of lines and irreducible conics on~$X$.
Here and below, we use the notation
\roster*
\item
$\Fn_1X$ for the graph of lines on~$X$,
\item
$\Fn_2X$ for the graph of all reduced conics on~$X$, and
\item
$\Fn_2\irr X\subset\Fn_2X$ for the induced subgraph of irreducible conics;
\endroster
in each graph, two vertices $u,v$ are connected by an edge of multiplicity
$u\cdot v$.
In addition to the Fano graphs~$\graph$ and
connected components of the respective \emph{absolute strata}
$\famX(\graph)$ in the space~$\famB$ of all Barth--Bauer octics, we also list
the \emph{relative strata} \smash{$\famXdm(\graph,\dm)\to\famX(\graph)$}
consisting of
pairs $(X,\dm)$, where $X$ is a Barth--Bauer octic and $\dm$ is a
distinguished unordered collection of Kummer conics on~$X$.

{\2
\convention
To avoid ambiguity, we emphasize that we consider \emph{smooth} octics only,
\ie, the polarization~$h$ is assumed very ample. By a \emph{conic} we mean a
reduced algebraic curve $C\subset X$
of arithmetic genus~$0$ and projective degree
$c\cdot h=2$. Thus, $C^2=-2$ and,
since exceptional divisors are not allowed,
for a conic~$C$ there are but two possibilities:
\roster*
\item
$C$ is irreducible, \ie, it is a true planar conic in~$\Cp5$, or
\item
$C=C_1+C_2$ splits into a pair of distinct intersecting lines{\3, so that
one has}
$C_1^2=C_2^2=-2$ and $C_1\cdot h=C_2\cdot h=C_1\cdot C_2=1$
(in particular, $C$ is still a planar curve).
\endroster
\endconvention}

The principal results of the paper, \viz. the complete list of deformation
classes, are collected in Tables~\ref{tab.1}--\ref{tab.3.2} (see
Theorems~\ref{th.codim.1}, \ref{th.codim.2}, \ref{th.codim.3}), itemized
according to the codimension of the strata in the $3$-parameter
family~$\famB$. (Following~\cite{degt:4Kummer,degt:conics}, we count both
conics and lines, hence both irreducible and reducible
conics.)
Here, in the introduction, we outline a few qualitative
consequences of this classification.

\theorem[see \autoref{proof.main}]\label{th.main}
The maximal number of conics on a Barth--Bauer octic is $176$.
Up to projective transformation, there is a unique
Barth--Bauer octic~$X_{176}$ with
$176$ conics, which are all irreducible\rom; it is given by
\[*
z^2_0+z^2_3-\phi z^2_4+\phi z^2_5=
z^2_1-\phi z^2_3+z^2_4-\phi z^2_5=
z^2_2+\phi z^2_3-\phi z^2_4+z^2_5=0,
\]
where $\phi:=(1+\sqrt5)/2$ is the golden ratio
\rom(see~\cite{Bonnafe.Sarti}\rom).
{\2This octic has no lines.}
\endtheorem

Recall that a typical smooth octic $K3$-surface in~$\Cp5$ is a \emph{triquadric},
\ie, a regular complete intersection of three quadrics. However, the moduli
space contains a divisor of \emph{special octics}, requiring at least one
cubic defining equation. Equivalently (Saint-Donat~\cite{Saint-Donat}),
special are the octics admitting an elliptic pencil of projective degree~$3$.
We assert that Barth--Bauer octics are never special.

\theorem[see \autoref{proof.triquadric}]\label{th.triquadric}
Any Barth--Bauer octic is a triquadric.
\endtheorem

Next, we support the speculation of~\cite{degt:conics} that, although it is
easier (at least, using the approach suggested in~\cite{degt:conics}) to count all,
not only irreducible
conics, all conics on a
polarized $K3$-surface are irreducible whenever their number is large
enough.

\theorem[see Tables~\ref{tab.1}--\ref{tab.3.2}]\label{th.bounds}
Let $X\subset\Cp5$ be a Barth--Bauer octic. Then\rom:
\roster*
\item
the maximal number of lines on~$X$ is $28$
\rom(a single octic, see \maxreal2 in \autoref{tab.3}\rom)\rom;
\item
the maximal number of reducible conics is $48$ \rom(same octic as above\rom)\rom;
\item
if $\ls|\Fn_2X|>128$, then $X$ is a singular $K3$-surface,
\ie, $\rank\NS(X)=20$\rom;
\item
if $\ls|\Fn_2X|>128$, then $X$ has no lines
\rom(hence, no reducible conics\rom)\rom;
\item
if $\ls|\Fn_2\irr X|>104$, then $X$ has no lines
\rom(hence, no reducible conics\rom).
\done
\endroster
\endtheorem

Theorems~\ref{th.main} and~\ref{th.bounds} should extend to all smooth octic
$K3$-surfaces, but the precise bounds may differ. For example, the sharp
upper bounds on the numbers of lines and reducible conics are essentially
found in~\cite{degt:lines}.

\theorem[see \cite{degt:lines} and \autoref{proof.lines}]\label{th.lines}
The maximal number of lines on a smooth octic $K3$-surface in~$\Cp5$ is $36$,
whereas the maximal number of reducible conics is $112$.
\endtheorem

\autoref{th.bounds} and the findings of \cite{degt:4Kummer,degt:conics}
suggest the following conjecture.

\conjecture\label{conj.N*}
There is a number $\Nstar{2n}(2)<\Nmax{2n}(2)$ with the following
property: if a smooth $2n$-polarized $K3$-surface $X\subset\Cp{n+1}$ has more
than $\Nstar{2n}(2)$ conics, then $X$ has no lines and, in particular, all
conics on~$X$ are irreducible.
\endconjecture

In conclusion, we address the question about the number of real conics on a real
surface (for which, as explained in~\cite{degt:4Kummer},
Barth--Bauer octics are not likely to provide good
examples). The current upper bound is as follows.

\theorem[see \autoref{proof.real}]\label{th.real}
The maximal number of real conics on a real Barth--Bauer octic
is $128$. There
is a unique $1$-parameter family of real Barth--Bauer octics with $128$ real
conics, see \maxreal1 in \autoref{tab.3}.
\endtheorem

\subsection{Digression: Mukai surfaces}\label{s.Mukai}
The second largest number of conics is~$160$ and, like $X_{176}$
in \autoref{th.main}, the
corresponding octic $X_{160}$ is also characterized by the presence of a
faithful projective symplectic action of a Mukai group~\cite{Mukai}, \viz.
$T_{192}$, see \autoref{ex.160}.
It is remarkable that \emph{Mukai surfaces} (\ie, $K3$-surfaces admitting a
faithful symplectic action of one of the eleven maximal groups in
\cite{Mukai}) maximize (sometimes conjecturally) the line or conic counts in
many degrees. For this reason, we use~\eqref{eq.Fn} and the known generic
N\'{e}ron--Severi lattices (see, \eg, \cite{Hashimoto}) to compute the Fano
graphs of all Mukai surfaces of degree $h^2\le8$. Results are shown in
\autoref{tab.Mukai}, where we list
\table
\caption{Conics on Mukai surfaces\noaux{ (see \autoref{s.Mukai})}}\label{tab.Mukai}
\def\*{\rlap{\smash{$^*$}}}
\def\?{\rlap{\smash{$^?$}}}
\centerline{\small\vbox{\halign{\strut\quad\hss$#$\hss\quad
 &\hss$#$\hss\quad&\hss$#$\hss\quad&\hss$#$\hss\quad&\hss$#$\hss\quad&#\hss\quad\cr
\noalign{\hrule\kern3pt}
G&(h^2,d)&\text{lines}&\text{conics}&T&Remarks\cr
\noalign{\kern1pt\hrule\kern3pt}
L_2(7)&(2,2)&&8526&[ 14, 0, 28 ]\cr
&(4,2)&&728&[ 14, 0, 14 ]\cr
\AG6&(2,2)&&8910\?&[ 12, 0, 30 ]&$=\Nmax2(2)$ ??\cr
&(6,2)&&285\*&[ 6, 0, 20 ]&$X_{285}$ in \cite{degt:conics}\cr
\SG5&(6,2)&&237&[ 10, 0, 20 ]&Table 7 in \cite{degt:conics}\cr
M_{20}&(4,2)&&800\?&[ 4, 0, 40 ]&Thm.\,1.3 in \cite{degt:4Kummer}\cr
&(8,2)&&176\?&[ 8, 4, 12 ]&\autoref{th.main}\cr
F_{384}&(4,1)&48&336+320&[ 8, 0, 8 ]\cr
&(8,1)&32&96+48&[ 4, 0, 8 ]\cr
\AG{4,4}&(8,2)&&144&[ 12, 0, 12 ]\cr
T_{192}&(4,1)&64\*&576+144&[ 8, 4, 8 ]&\autoref{rem.X64}\cr
&(8,2)&&160&[ 4, 0, 24 ]&\autoref{ex.160}\cr
H_{192}&(4,1)&48&336+168&[ 8, 0, 12 ]\cr
&(8,1)&32&96+12&[ 4, 0, 12 ]\cr
N_{72}&(6,2)&&225&[ 6, 0, 36 ]\cr
M_9&(2,1)&144\*&5112+2988&[ 12, 6, 12 ]&$=\Nstar2(2)$ ??\cr
T_{48}&(2,1)&108&2862+3180&[ 16, 8, 16 ]\cr
\noalign{\kern1pt\hrule}
\crcr}}}
\endtable
\roster*
\item
the Mukai group~$G$ (in the notation of~\cite{Mukai}),
the degree~$h^2$ of the
model,
and its \emph{depth}
{\3$d$ defined \via\ }
$d:=\gcd\{x\cdot h\,|\,x\in\NS(X)\}$,
\item
the numbers of lines and conics on~$X$; the latter is shown as a single count
if all conics are irreducible, or as
$(\text{reducible})+(\text{irreducible})$ otherwise,
\item
the transcendental lattice $T(X)$.
\endroster
We omit hyperelliptic models (except degree $h^2=2$) and those of
depth $d>2$ (as they obviously have no lines or conics). The line/conic
counts known or conjectured to be maximal are marked with $^*$ or
$^?$, respectively.

Some of these configurations have already appeared elsewhere (see the remark
column),
whereas others seem to be new. Probably, the most important discovery is the
following observation (see $\AG6$ and $M_9$ in \autoref{tab.Mukai};
\cf. \autoref{conj.N*}).

\observation
One has $\Nmax2(2)\ge8910$ and $\Nstar2(2)\ge8100$ (if defined).
\endobservation

\subsection{Contents of the paper}
In \autoref{S.BB}, we recall a few basic facts about (polarized) Kummer
surfaces (\autoref{s.prelim}, \autoref{s.BB}),
analyze a very general Barth--Bauer octic (\autoref{s.very.general}), and lay
the basis for the study of the equiconical strata of positive
codimension (\autoref{s.components}, \autoref{s.supports}). At the end, in
\autoref{proof.triquadric}, we prove \autoref{th.triquadric}.

In \autoref{S.codim.1} we perform a deep case-by-case analysis resulting in
the five codimension~$1$ strata listed in \autoref{tab.1}, see
\autoref{th.codim.1}.
Finally, in \autoref{S.codim.large}, we list all strata of higher codimension
(see Theorems~\ref{th.codim.2}, \ref{th.codim.3} and
Tables~\ref{tab.2}--\ref{tab.3.2}) and give formal proofs of the principal
results of the paper stated in the introduction.

\subsection{Acknowledgements}
This paper was finalized during my sabbatical stay at the
\emph{Max-Planck-Institut f\"{u}r Mathematik}, Bonn; I am grateful to this
institution for its support and the excellent working environment.

\section{Barth--Bauer octics}\label{S.BB}

In this section, we recall a few basic facts about (polarized) Kummer
surfaces (see \autoref{s.prelim} and \autoref{s.BB}),
analyze a very general Barth--Bauer octic (see \autoref{s.very.general}), and lay
the basis for the study of the equiconical strata of positive
codimension (see \autoref{s.components} and \autoref{s.supports}). In
\autoref{proof.triquadric}, we use the machinery of \autoref{s.supports} to
prove \autoref{th.triquadric}.

\subsection{Preliminaries}\label{s.prelim}
Let $\dm$ be a $16$-element set; denote $\CC_0:=\{\varnothing\}$,
$\CC_{16}:=\{\dm\}$.
A \emph{Kummer structure} on~$\dm$ is a
collection $\CO_8$ of $30$ eight-element subsets $\go\subset\dm$ such that
$\CO_*:=\CO_0\cup\CO_8\cup\CO_{16}$ is closed under the
symmetric difference~$\sd$. (Here and below, for a subset~$\CS_*$ of a
power set, we use the convention $\CS_n:=\{\go\in\CS_*\,|\,\ls|\go|=n\}$.
According to Nikulin~\cite{Nikulin:Kummer}, any Kummer structure is
standard: there is a bijection between~$\dm$
and a codeword of length~$16$ of the
(extended binary) Golay code~$\Golay_*$
(see, \eg, \cite{Conway.Sloane}) such
that $\CO_8=\{\go\in\Golay_8\,|\,\go\subset\dm\}$.
Then one also has
\[*
\CC_*:=\bigl\{\gs\subset\dm\bigm|
   \text{$\ls|\gs\cap\go|=0\bmod2$ for all $\go\in\CO_*$}\bigr\}
 =\bigl\{\gs\cap\dm\bigm|\gs\in\Golay_*\bigr\},
\]
and the setwise stabilizer of~$\CO_*$ in $\SG{16}$ is the restriction
to~$\dm$ of its stabilizer in the Mathieu group~$M_{24}$.
This group acts transitively on~$\CC_4$ and, hence, on the set of
\emph{$8$-Kummer structures} (\cf. \cite{degt:4Kummer}) defined \via
\[*
\CK_*:=\CK_*(\gk):=\bigl\{\gk\sd\go\bigm|\go\in\CO_*\bigr\}\quad
 \text{for some fixed $\gk\in\CC_4$}.
\]
Note that $\CK_*$ is generated by any of the four elements $\gk\in\CK_4$ and
$\CO_*$ is recovered back from~$\CK_*$ \via\
$\CO_*=\{\gr\sd\gs\,|\,\gr,\gs\in\CK_*\}$.
The setwise stabilizer~$\G$ of $\CK_*$ is a group of order $9216$.

Throughout the paper, we use the following shortcuts (where
$\gr,\gs\subset\dm$):
\[*
\hbar:=\tfrac12h\in\Q h,\qquad
\gs:=\sum_{e\in\gs}e\in\Z\dm,\qquad
\ds\gs\gr:=\tfrac12(\gs\cap\gr)-\tfrac12(\gs\sminus\gr)\in\Q\dm.
\]
The other terminology and notation related to lattices is quite standard,
\cf.~\cite{degt:4Kummer}.

From now on, we fix an $8$-Kummer structure~$\CK_*$ and consider the lattices
\begin{alignat}2
\label{eq.L}
&\L         &&:=2\bE_8\oplus3\bU\cong H_2(X;\Z)\
 \text{for a $K3$-surface~$X$};\\
\label{eq.S}
&\bS        &&:=\bS(\CO_*)\supset\Z\dm\
 \text{is the extension \via\ all $\ds\go\varnothing$, $\go\in\CO_*$};\\
\label{eq.T}
&\Sperp     &&:=\bS_\L^\perp\cong3\bU(2)\
 \text{for a fixed primitive isometry $\bS\into\L$};\\
\label{eq.Sh}
&\lattice   &&:=\lattice(\CK_*)\supset\bS+\Z h\
 \text{is the extension \via\ all $\hbar+\ds\gk\varnothing$, $\gk\in\CK_*$}.
\end{alignat}
A primitive isometry $\bS\into\L$ in~\eqref{eq.T} is unique up to isomorphism
(see~\cite{Nikulin:Kummer}), and in
\eqref{eq.Sh} we let $h\cdot e=2$ for $e\in\dm$.
In particular, \eqref{eq.T} implies that
\[
\text{$u^2=0\bmod4$,\quad $u\cdot v=0\bmod2$\quad for any $u,v\in\Sperp$}.
\label{eq.Sperp}
\]
We also introduce the equivalence relations
\[*
\text{$\gr\eq\gs$ iff $\gr\sd\gs\in\KC$},\qquad
\text{$\gr\Eq\gs$ iff $\gr\sd\gs\in\KC\cup\KK$}
\]
on~$\CS_*$
and respective equivalence classes $\cl\cdot$, $\Cl\cdot$.

The \emph{parity}
of a set $\gs\in\CC_*$ is 
$\ls|\gs\cap\gk|\bmod2$
for some (equivalently, any) $\gk\in\CK_*$.
Since any $\gs\in\CC_*\cup\CK_*$ is even,
the parity is preserved by~$\eq$ and~$\Eq$.
The $\G$-action on~$\CC_*$ respects $\bar\ $,
$\sd$, parity, and both $\eq$ and $\Eq$; its orbits are shown in
\autoref{tab.orbits}, where most nonempty cells represent a single orbit each,
shown as
$\#(\text{$\eq$-classes})\times\ls|\text{class}_n|$.
The two
exceptional cases $(\operatorname{even})_8$ and $\CO_8$ consist of two
$\G$-orbits each: an extra invariant of a set~$\gs$ is the existence of
$\gk\in\CK_4$ such that $\gk\subset\gs$.
However,
the induced actions on $(\operatorname{even})_8/{\eq}$ and $\CO_8/{\eq}$ are
still transitive.
\table
\caption{$\G$-orbits on $\CC_n$}\label{tab.orbits}
\def\0{\phantom0}
\centerline{\small\vbox{\halign{&\strut\quad\hss$#$\hss\quad\cr
\noalign{\hrule\kern3pt}
n&\text{even}&\text{odd}&\KC&\KK\cr
\noalign{\kern1pt\hrule\kern3pt}
0&&&1\times1\cr
4&18\times4\0&16\times4\0&&1\times4\0\cr
6&12\times16&16\times16\cr
8&18\times(8+16)&16\times24&1\times(6+24)&1\times24\cr
10&12\times16&16\times16\cr
12&18\times4\0&16\times4\0&&1\times4\0\cr
16&&&1\times1\cr
\noalign{\kern1pt\hrule}
\crcr}}}
\endtable

The next lemma is a straightforward application
of~\cite{Nikulin:Kummer,Nikulin:forms}.
We present a partial statement which is used in this paper; more
details are found in~\cite{degt:4Kummer}.

\lemma\label{lem.kernels}
For a Kummer structure $\CO_*$ and primitive isometry
$\bS:=\bS(\CO_*)\into\L$, consider an overlattice $\bS\subset N\subset\L$
primitive in~$\L$ and
let $\bS^\perp:=N\cap\Sperp$.
Then, for each vector $u\in\bS^\perp$, there is a class $\Cal U\in\dm/{\eq}$
such that, for each $\gu\in\Cal U$,
\[*
2\ls|\gu|=u^2\bmod8\quad\text{and}\quad\tfrac12(u+\gu)\in N.
\def\qedsymbol{\donesymbol}\pushQED\qed\qedhere
\]
\endlemma

\subsection{Barth--Bauer surfaces}\label{s.BB}
According to Nikulin~\cite{Nikulin:Kummer},
a Kummer surface $(X,\dm)$ defines a canonical Kummer structure on the set~$\dm$ of
its Kummer divisors, and
the N\'{e}ron--Severi lattice $\NS(X)$ is a primitive extension
of~$\bS$ in~\eqref{eq.S}. If $X$ is polarized, $\NS(X)\ni h$, so that each
Kummer divisor $e\in\dm$ is a conic, $e\cdot h=2$, then
\[
\L\supset\NS(X)\supset\bS+\Z h=\bS\oplus\Z\th,\qquad
 \th:=\dm+h,\quad\th^2=32+h^2;
\label{eq.th}
\]
in particular, $h^2=0\bmod4$ by~\eqref{eq.Sperp}.

From now on, we assume that $h$ is very ample and $h^2=8$, even though some
formulas below are written for arbitrary~$h^2$.
By Saint-Donat~~\cite{Saint-Donat}, neither $h$ nor~$\th$ is divisible
by~$2$ in $\NS(X)$;
hence, the class $\Cal U\in\dm/{\eq}$ given by
\autoref{lem.kernels} for $u=\th$ is a certain $8$-Kummer structure~$\CK_*$,
so that $\NS(X)$ is a primitive extension of the lattice $\lattice$
in~\eqref{eq.Sh}.
This extension must be \emph{geometric} in the following sense.

\definition[\cf. Saint-Donat~\cite{Saint-Donat}]\label{def.admissible}
A hyperbolic overlattice
$N\supset\Z\dm+\Z h$
is called
\emph{admissible} if
\roster
\item\label{i.Veronese}
$h$ is not divisible by~$2$ in~$N$, and
\endroster
there is no vector $r\in N$ such that either
\roster[\lastitem]
\item\label{i.exceptional}
$r^2=-2$ and $r\cdot h=0$ (\emph{exceptional divisor}), or
\item\label{i.isotropic}
$r^2=0$ and $r\cdot h=\pm2$ (\emph{$2$-isotropic vector}), or
\item\label{i.missing}
$r^2=-2$, $r\cdot h=1$, and $r\cdot e<0$ for some $e\in\dm$
(\emph{missing conic}).
\endroster
An admissible lattice~$N$ is called \emph{geometric} if
the isometry $\bS\into\L$, see~\eqref{eq.T},
extends to a primitive
isometry $N\into\L$.
\enddefinition

\remark\label{rem.admissible}
According to Nikulin~\cite{Nikulin:Kummer}, in any geometric overlattice
$N\supset\Z\dm+\Z h$ one
{\3necessaroly}
has $N\cap(\Q\dm+\Q h)=\lattice(\KK)$ for some
$8$-Kummer structure~$\KK$ on~$\dm$. For this reason we usually fix~$\KK$ and
work with overlattices of~$\lattice$.
\endremark

Conversely, a standard chain of arguments based on the global Torelli
theorem \cite{Pjatecki-Shapiro.Shafarevich}, surjectivity of the period
map~\cite{Kulikov:periods},
and the results of Nikulin~\cite{Nikulin:finite.groups} and
Saint-Donat \cite{Saint-Donat} shows that each geometric overlattice
$N\supset\lattice\supset\dm$ serves as $\NS(X)$ for some Barth--Bauer
octic $(X,\dm)$. Indeed, an
abstract $K3$-surface~$X$ is given by the surjectivity of the period map;
then,
conditions~\iref{i.isotropic} and~\iref{i.Veronese} assert that the linear
system~$h$ defines a map
$\Gf_h\:X\to\Cp5$ which is
birational onto its image, condition~\iref{i.exceptional} makes the image
$\Gf_h(X)$ smooth, and condition~\iref{i.missing} is equivalent to
the requirement that each class $e\in\dm$ {\2should} represent an
\emph{irreducible} $(-2)$-curve on~$X$.

The moduli space of octics~$X$ obtained in this way is discussed in
\autoref{s.components} below. The Fano graphs of~$X$ (see \autoref{S.intro})
can be computed in terms of the polarized lattice $N:=\NS(X)\ni h$ using the
description of the nef cone in Huybrechts \cite[\S\,8.1]{Huybrechts} and
Vinberg's algorithm~\cite{Vinberg:polyhedron} (\cf. also
\cite{degt:4Kummer,degt.Rams:octics}): identifying $(-2)$ curves on~$X$ with
their classes in~$N$, we have
\[
\aligned
\Fn_n(N,h)&:=\bigl\{u\in N\bigm|
 \text{$u^2=-2$ and $u\cdot h=n$}\bigr\},\quad n=1,2,\\
\Fn_2\irr(N,h)&:=\bigl\{u\in\Fn_2(N,h)\bigm|
 \text{$u\cdot v\ge0$ for all $v\in\Fn_1(N,h)$}\bigr\}.
\endaligned
\label{eq.Fn}
\]
The inverse of~\eqref{eq.Fn} assigns to a bi-colored graph~$\graph$ the
$8$-polarized lattice
\[
\Fano(\graph):=(\Z\graph+\Z h)/{\ker},\quad
 h^2=8,\quad
 \text{$h\cdot v=\operatorname{color}(v)$ for $v\in\graph$},
\label{eq.Fano}
\]
where $\Z\graph$ is freely generated by the vertices $v\in\graph$ and
$u\cdot v=n$ whenever $u,v\in\graph$ are connected by an $n$-fold edge.
{\2(Here, $\ker(\cdot):=(\cdot)^\perp$ refers to the kernel of the bilinear
form.)}
\latin{A priori}, $\Fano(\graph)$ is neither geometric nor admissible;
in fact, it does not even need to be hyperbolic.

\subsection{Generic Barth--Bauer octics}\label{s.very.general}
A very general Barth--Bauer
octic $X\subset\Cp5$ has the minimal N\'{e}ron--Severi lattice
$\NS(X)=\lattice$,
and a computation using~\eqref{eq.Fn} shows that $X$ has exactly $32$ conics,
all irreducible:
\roster*
\item
the $16$ original Kummer conics $e\in\dm$, and
\item
$16$ pairwise disjoint irreducible \emph{Barth--Bauer}, or \emph{\BB-conics}
\[
\hh+\ds\gk\gs,\qquad \gs\subset\gk\in\CK_4,\ \ls|\gs|=1;
\label{eq.BB}
\]
these conics have pattern $\cn12-3$ in the notation of
\eqref{eq.pattern} below.
\endroster

\remark\label{rem.index=4}
Note that $\lattice=\NS(X)$
{\3(see the first row in \autoref{tab.1})}
is not generated over~$\Z$ by $h$ and conics: one
has $[\lattice:\Fano(\Fn\lattice)]=4$.
(There are but two other strata with this property, see Tables~\ref{tab.1}
and~\ref{tab.2}.)
It is for this reason that the group $\OG_h(\lattice)=\G\times\Z/2$ is much
smaller than the full group $\Aut(\Fn\lattice)$.
\endremark

\table
\caption{Symplectic groups $G_\Go$}\label{tab.groups}
\def\0{\phantom0}%
\makeatletter
\centerline{\small\offinterlineskip\vbox{%
\def\\#1(#2,#3){\omit\raise6pt\hbox{\hypertarget{sym(#2,#3)}{}}\sep\hss$#1$\hss\sep
 \protected@write\@auxout{}{\string\newcs{(#2,#3)}{\string\hyperlink{sym(#2,#3)}{#1}}}%
 &#2&#3}%
\let\sep\quad
\halign{&\strut\sep\hss$#$\hss\sep\cr
\noalign{\hrule\kern3pt}
\#&\ls|G_\Go|&\text{index}&G_\Go\cr
\noalign{\kern1pt\hrule\kern3pt}
\\21(16,14)&C_4^2\cr
\\39(32,27)&2^4C_2\cr
\\49(48,50)&2^4C_3\cr
\\77(192,1493)&T_{192}\cr
\\81(960,11357)&M_{20}\cr
\noalign{\kern1pt\hrule}
\crcr}}}
\endtable
Denote by $\G_\Go\cong(\Z/4)^4$
(see \#\symplectic(16,14) in \autoref{tab.groups}) the subgroup
of~$\G$ acting identically on~$\dm/{\eq}$. Clearly,
$\G_\Go\subset\OG_h(\lattice)$ is the subgroup acting identically on
$\discr\lattice$; this action extends to any overlattice $N\supset\lattice$
and, by the global Torelli theorem, gives rise to a projective symplectic
action on any Barth--Bauer octics. All extensions of $\G_\Go$ acting
symplectically on (generic in their respective strata) Barth--Bauer octics
are listed in \autoref{tab.groups}, where \#\ and ``index'' refer,
respectively, to the list in Xiao~\cite{Xiao:Galois} and \GAP~\cite{GAP4}
small group library and the last column is the notation
in~\cite{Xiao:Galois}.

\subsection{Connected components}\label{s.components}
Given a lattice~$L$,
we denote by $\OGplus(L)$ the group of auto-isometries of~$L$
preserving a \emph{positive sign structure}, \ie, coherent
orientation of all maximal positive definite subspaces of $L\otimes\R$.

Let $N\supset\lattice$ be a geometric overlattice, see \autoref{def.admissible},
and $G\subset\OG_h(N)$ a fixed subgroup: in what follows, we will have either
$G=\OG_h(N)$ or $G=\stab\dm$.
Two isometries $\Gf_i\:N\into\L$, $i=1,2$, are said to be
\emph{$G$-equivalent} if there exists a pair of isometries
$g\in G$, $f\in\OGplus(\L)$ such that $f\circ\Gf_1=\Gf_2\circ g$.

Fix a bi-colored graph~$\graph$ and consider geometric finite index extensions
\[
N\supset\Fano(\graph)\ni h\quad\text{such that}\quad\Fn(N,h)=\graph.
\label{eq.N}
\]
Using Dolgachev's~\cite{Dolgachev:polarized} coarse moduli space of lattice
polarized $K3$-surfaces and factoring out the projective group, one easily
concludes (see~\cite{degt:conics}) that the connected components of the equiconical stratum
$\famX(\graph)$ are of the form $\famX(N\into\L)$, where
\roster*
\item
$N\supset\Fano(\graph)\ni h$ is a geometric finite index extension as
in~\eqref{eq.N}, regarded up to lattice isomorphism preserving~$h$, and
\item
$N\into\L$ is an $\OG_h(N)$-equivalence class of primitive isometries.
\endroster
A similar statement holds for the relative stratum $\famXdm(\graph,\dm)$,
except that
\roster*
\item
$N$ is regarded up to isomorphism preserving~$h$ and $\dm$ (as a set), and
\item
$N\into\L$ is a $(\stab\dm)$-equivalence class of primitive isometries.
\endroster
In both cases, a component $\famX(\Gf\:N\into\L)$ is real if and only if
$\Gf$ is equivalent to $g\circ\Gf$ for some (equivalently, any)
$g\in\OG(\L)\sminus\OGplus(\L)$.

Thus, the connected components of the strata associated to a graph~$\graph$
are in a bijection with the
equivalence
classes of the diagrams
\[
\Fano(\graph)\into N\into\L,
\label{eq.diag}
\]
where $N$ is admissible, the former
arrow is a finite index extension as in~\eqref{eq.N}, and the latter arrow is a
primitive isometry. {\2The equivalence is up to the group $\OG^+(\L)$ of
auto-isometries preserving the orientation of maximal positive definite
subspaces and appropriate, depending on the kind of strata considered,
subgroup of the group $\OG_h(N)$ of autoisometries of~$N$ preserving the
polarization~$h$.
}

\subsection{The computation\noaux{ (see Nikulin~\cite{Nikulin:forms})}}\label{s.computation}
At each step\mnote{\2new section}
in~\eqref{eq.diag}, there are but a finite number of
choices, easily found in terms of the \emph{discriminant group}
\[*
\discr S:=S\dual\!/S,\quad\text{where}\quad
 S\dual=\bigl\{y\in S\otimes\Q\bigm|\text{$x\cdot y\in\Z$ for all $x\in S$}\bigr\},
\]
equipped with
the $\Q/2\Z$ quadratic form $(y\bmod S)\mapsto y^2\bmod2\Z$
(see~\cite{Nikulin:forms},
where the less distinctive notation $q_S$ is used).
This group is abelian and finite.

Given a nondegenerate even lattice~$S$, the map
$N\mapsto\CK:=N\bmod S$
establishes a bijection between the isomorphism classes of
finite index extensions $N\supset S$ and
isotropic subgroups $\CK\subset\discr S$.
Furthermore,
\[
\text{$g\in\OG(S)$ extends to~$N$ if and only if $g(\CK)=\CK$}.
\label{eq.g.ext}
\]
We mainly work with polarized hyperbolic lattices $S\ni h$ rationally
generated by lines and conics. In this case we obviously have
$\OG_h(S)\subset\Aut\graph=\OG_h\bigl(\Fano(\graph)\bigr)$,
where $\graph:=\Fn(S,h)$; the latter group is
computed using \texttt{Digraphs} package in \GAP~\cite{GAP4},
and the former
is given by~\eqref{eq.g.ext}.

Thus, we can effectively list the isomorphism classes of the finite index
extensions $N\supset S\ni h$. For each polarized lattice $N\ni h$
obtained, we check the admissibility (see \autoref{def.admissible})
and make sure that $\Fn(N,h)=\Fn(S,h)$, as otherwise we would have started
from the larger graph $\Fn(N,h)$ in the first place.


The second arrow in~\eqref{eq.diag} is a primitive extension. Given such an
extension $N\into\L$, the genus $g(T)$ of the \emph{transcendental lattice}
$T:=N^\perp$ is determined by $\discr T\cong-\discr N$. (In particular, the
existence of a primitive extension depends on $\rank N$ and $\discr N$ only, see
\cite[Theorem 1.12.2]{Nikulin:forms}.) The extensions with a given
lattice~$T$ are in a bijection with anti-isometries
$\Gf\:\discr N\to\discr T$; furthermore
\[
\text{$g\in\OG(N)$ and $h\in\OG(T)$ extend to~$\L$ if and only if
 $\Gf\circ g=h\circ\Gf$}.
\label{eq.gh.ext}
\]
Thus, to list the
isomorphism classes of primitive extensions $N\into\L$, we need to
\roster
\item\label{i.T.genus}
list (representatives of) the isomorphism classes $T\in g(T)$, and
\item\label{i.T.cosets}
for each class~$T$, compute the quotient
$\OG_h(N)\backslash\!\Aut(\discr T)/\!\OG^+(T)$, which makes sense upon
fixing an anti-isometry $\discr N\to\discr T$.
\endroster
Considering that $\OG_h(N)$ is known, for~\iref{i.T.cosets} we
merely need the image
of the canonical homomorphism $\OG^+(T)\to\Aut(\discr T)$.

If $\rank T=2$ and, hence, $T$ is positive definite, we use Gauss'
theory~\cite{Gauss:Disquisitiones} of binary quadratic forms (see also
\cite[Theorem 56]{Dickson}): the reduced form suggested therein lets one list
all representatives of a genus and compute the finite groups $\OG^+(T)$.

If $\rank T\ge3$, we use Miranda--Morrison theory~\cite{Miranda.Morrison:book},
which combines both
{\3the genus group} $g(T)$ and
{\3cokernel}  $\Coker\bigl[\OG^+(T)\to\Aut(\discr T)\bigr]$
in a single abelian group $E(T)$ that is computed in terms of the
discriminant $\discr T$. A brief account of the theory is found
in~\cite{degt:geography}; for the computation details, we refer to
\cite[Chapter VII]{Miranda.Morrison:book}.


\subsection{The supports of a vector}\label{s.supports}
In view of \autoref{rem.admissible},
the graphs~$\graph$ to be tried for~\eqref{eq.diag} are of the form
$\graph:=\Fn\lattice\(u_i)$, where $\lattice\(u_i)\supset\lattice$ is a
primitive corank~$r$ extension generated by $r$ extra lines or conics
$u_1,\ldots,u_r$.
The next lemma {\2(some parts of which are obvious geometrically)}
controlls such extensions by bounding
the intersection indices of lines and conics.

\lemma\label{lem.intersection}
Let $X\subset\Cp5$ be a smooth $K3$-octic, $l_1,l_2\in\NS(X)$ a pair of
distinct lines
on~$X$, and $c_1,c_2\in\NS(X)$ a pair of distinct conics. Then one has
\[*
l_1\cdot l_2\le1,\quad
l_1\cdot c_1\le2\ (\text{or $1$, if $X$ is a triquadric}),\quad
c_1\cdot c_2\le2.
\]
\endlemma

\proof
By the Hodge index theorem, the lattice $\NS(X)$ is hyperbolic. Hence, for
any pair of vectors $u,v\in\NS(X)$, one has
\[
\det(\Z h+\Z u+\Z v)\ge0,
\label{eq.det}
\]
with the equality attained if and only if $h,u,v$ are linearly dependent.

Applying~\eqref{eq.det} to one of the three pairs in the statement, we obtain
\[*
l_1\cdot l_2\le2,\qquad
l_1\cdot c_1\le2,\qquad
c_1\cdot c_2\le3,
\]
and there remains to rule out the possibilities $l_1\cdot l_2=2$ and
$c_1\cdot c_2=3$.

In the former case, $l_1\cdot l_2=2$, the lattice contains the $2$-isotropic
vector  $l_1+l_2$, see \autoref{def.admissible}\iref{i.isotropic}, and the
map $X\to\Cp5$ defined by~$h$ is two-to-one, see~\cite{Saint-Donat}.

In the latter case, $c_1\cdot c_2=3$, the determinant~\eqref{eq.det} vanishes
and we obtain a relation $2h=c_1+c_2$. Hence, $h$ is divisible by~$2$ in
$\NS(X)$ and the map $X\to\Cp5$ is also two-to-one, factoring through the
Veronese embedding $\Cp2\into\Cp5$,
see~\cite{Saint-Donat}.

For the bound $l_1\cdot c_1\le1$, observe that, if $l_1\cdot c_1=2$, then the vector
$e:=l_1+c_1$ is \emph{$3$-isotropic}: $e^2=0$, $e\cdot h=3$. According
to~\cite{Saint-Donat} (see also \cite{degt.Rams:octics}),
the presence of such a vector in $\NS(X)$ is
equivalent to the fact that $X$ is special.
\endproof

\table
\caption{Sylvester test for conics (left) and lines (right)}\label{tab.Sylvester}
\def\bx#1{\hbox to9pt{\hss$#1$\hss}}%
\def\select#1{\if#1.\cdot\else
 \ifcase#1
  \times\or
  \bullet\or
  \circ\or
  \circ\fi\fi
}%
\def\r#1{\hss$\scriptstyle#1$\ \,}%
\def\q#1{\omit\bx{\scriptstyle#1}}%
\hrule height0pt
\centerline{\baselineskip10pt%
\vtop{\halign{#&&\bx{\select{#}}\cr
&\q{0}&&\q{2}&&\q{4}&&\q{6}&&\q{8}&&\q{10}&&\q{12}&&\q{14}&&\q{16}\cr
\r{0} &2&2&.&.&.&.&.&.&.&.&.&.&.&.&.&.&.\cr
\r{2} &0&.&.&.&.&.&.&.&.&.&.&.&.&.&.\cr
\r{4} &1&.&.&.&.&.&.&.&.&.&.&.&.\cr
\r{6} &1&.&.&.&.&.&.&.&.&.&.\cr
\r{8} &1&.&.&.&.&.&.&.&.\cr
\r{10}&1&.&.&.&.&.&.\cr
\r{12}&1&2&2&.&2\cr
\r{14}&0&0&0\cr
\r{16}&2\cr}}%
\quad\
\vtop{\halign{#&&\bx{\select{#}}\cr
&\q{0}&&\q{2}&&\q{4}&&\q{6}&&\q{8}&&\q{10}&&\q{12}&&\q{14}&&\q{16}\cr
&0&0&.&.&.&.&.&.&.&.&.&.&.&.&.&.&.\cr
&0&.&.&.&.&.&.&.&.&.&.&.&.&.&.\cr
&1&.&.&.&.&.&.&.&.&.&.&.&.\cr
&1&.&.&.&.&.&.&.&.&.&.\cr
&1&.&.&.&.&.&.&.&.\cr
&1&.&.&.&.&.&.\cr
&1&.&.&.&.\cr
&0&0&0\cr
&0\cr}}}
\endtable

In view of \autoref{lem.intersection}, if $e$ is an irreducible
conic on~$X$, then $u\cdot e\in\{0,1,2\}$ for any line or conic $u\ne e$. It
follows that a $1$-vector extension
\[*
\lattice\(u):=(\lattice+\Z u)/{\ker}
\]
(not necessarily proper) is uniquely determined by the degree $u\cdot h$
and two \emph{supports}
\[*
\supp_iu:=\bigl\{e\in\dm\bigm|u\cdot e=i\bigr\}\subset\dm,\quad i=1,2,
\]
which are two disjoint subsets of~$\dm$. Letting $p:=\ls|\supp_1u|$ and
$q:=\ls|\supp_2u|$, we will say that
\[
\text{$u$ has \emph{pattern} $\ln p-q$ (if it is a line) or $\cn p-q$ (if it is a conic)}.
\label{eq.pattern}
\]
Assuming that $\lattice\(u)$ is an integral lattice, we also have
\[
\text{$\supp_1u\in\CC_*$ is an even (resp.\ odd) set if $u\cdot h$ is even (resp.\ odd)}.
\label{eq.parity}
\]
Finally, denoting by $u_\bS$ the orthogonal projection of~$u$
to~$\lattice\otimes\Q$, we find that
\[
u_\bS^2=-\frac{p}2-2q+\frac{(p+2q+\Ge)^2}{h^2+32},
\label{eq.proj}
\]
where $p,q$ are as above and $\Ge:=u\cdot h$.
The lattice $\lattice\(u)$ is
hyperbolic and of corank~$1$ over~$\lattice$ if and only if $u_\bS^2>u^2=-2$. This
inequality results in \autoref{tab.Sylvester}
(the pairs marked with a $\cdot$ are ruled out), where, in view
of~\eqref{eq.parity}, only even values of~$p$ are considered.
For the reader's convenience,
the pairs ruled
out by~\eqref{eq.parity} and \autoref{tab.orbits} are marked with a $\times$,
and those prohibited in
\autoref{s.restrictions} below
are marked with a $\circ$.

\subsection{Proof of \autoref{th.triquadric}}\label{proof.triquadric}
As already mentioned, \cite{Saint-Donat}
(see also \cite{degt.Rams:octics}) states that a smooth $K3$-octic is special if and only
if the lattice $\NS(X)$ contains a \emph{$3$-isotropic vector}, \ie, a
vector~$u$ such that $u^2=0$ and $u\cdot h=3$. Applying~\eqref{eq.det} to
$v=e\in\dm$,
we get $u\cdot e\in\{0,1,2\}$. Hence, similar to \autoref{s.supports},
an extension $\lattice\(u)$
by a $3$-isotropic vector~$u$ is determined by the pair of supports
$\supp_iu\subset\dm$, $i=1,2$.
Arguing as in \autoref{s.supports}, we arrive at $u_\bS^2\ge u^2=0$,
where $u_\bS^2$ is given by~\eqref{eq.proj} with $\Ge=3$.
This inequality results in $\ls|\supp_1u|\in\{0,14,16\}$.
On the other hand,
$\supp_1u\in\CC_*$ is an odd set, see~\eqref{eq.parity}, contradicting to
\autoref{tab.orbits}.
\qed

\section{Strata of codimension~$1$}\label{S.codim.1}

The goal of this section is the description of the codimension~$1$ strata in
the space~$\famB$ of Barth--Bauer octics. The following theorem is proved in
\autoref{proof.codim.1} below.

\theorem\label{th.codim.1}
The space $\famB$ has five irreducible equiconical strata of codimension~$1$,
\viz.\ those listed in \autoref{tab.1}.
Each stratum consists of a single real component.
\table
\caption{Strata of codimension $\le1$ \rm(see~\autoref{th.codim.1})}\label{tab.1}
\centerline{\vbox{\small\rm\offinterlineskip
\tabdefs
\def\frac#1#2{#1/#2}%
\tabskip0pt plus 1fil
\halign{\strut\sep\hss$#$\hss\sep\tabskip0pt&
 \strut\sep\hss$#$\hss\sep&
 \strut\sep\hss$#$\hss\sep&
 \strut\sep\hss$#$\hss\sep&
 \strut\sep\hss$#$\hss\sep&
 \strut\sep\hss$#$\hss\sep&
 \strut\sep\hss$#$\hss\sep&
 \strut\sep\hss$#$\hss\sep&
 \strut\sep\hss$#$\hss\sep&
 \strut\sep\hss$#$\hss\sep&
 \strut\sep\hss$#$\hss\quad\tabskip0pt plus 1fil\cr
\noalign{\hrule\kern3pt}
\text{Name}&\text{Patterns}&\Gd_2^2&\Gd_5&\text{Lines}&\text{Conics}&\ls|G|&
 i_\dm&G_\Go&\ls|{\det}|&(r,c)\cr
\noalign{\kern1pt\hrule\kern3pt}
\text{open}&&&&&32&18432\cdot864&2&\(16,14)\^1&640\^4&(1,0)\EQ(1,0)\cr
\noalign{\kern1pt\hrule\kern3pt}
\astrat1&\ln4-0&\frac58&0&4&32&1152&2&\(16,14)\^1&400&(1,0)\EQ(1,0)\cr
\astrat2&\ln6-0,\ln12-0,\cn4-0&\frac58&\pm1&20&16\+20&576&1&\(16,14)\^1&144&(1,0)\EQ(1,0)\cr
\noalign{\kern1pt\hrule\kern3pt}
\astrat3&\cn4-0,\cn12-0&\frac24&\pm3& &40&1024\cdot16&2&\(16,14)\^1&576\^2&(1,0)\EQ(1,0)\cr
\astrat4&\cn6-0,\cn10-0&\frac22&\pm4& &64&3072&4&\(16,14)\^1&384&(1,0)\EQ(1,0)\cr
\astrat5&\cn8-0&\frac24&0& &80&2048&2&\(16,14)\^2&320&(1,0)\EQ(1,0)\cr
\noalign{\kern1pt\hrule}
\crcr}}}
\endtable
\endtheorem

For completeness, in the first row of
\autoref{tab.1} we also show the open stratum
of codimension~$0$, \ie, the one consisting of generic
Barth--Bauer octics.

\subsection{Notation in Tables
\ref{tab.1}--\ref{tab.3.2}}\label{s.list.1}
The rows of each table represent the isomorphism classes of pairs
$(\graph,\dm)$, where $\graph$ is a Fano graph and $\dm\subset\graph$ is a
distinguished set of $16$ irreducible Kummer conics.
The rows corresponding to isomorphic abstract bi-colored
graphs~$\graph$ are prefixed with equal
superscripts.
Listed in \autoref{tab.1} are
\roster*
\item
the name of the stratum (for further references),
\item
the patterns of the extra lines and conics, see~\eqref{eq.pattern}, and
\item
a description of the images $\Gd_p(u)\in\discr_p\lattice$, $p=2,5$ (see
\autoref{s.clusters} below), of a distinguished generator~$u$.
\endroster
Instead, the first column of the other tables merely lists
\roster*
\item
the types of the clusters (see \autoref{s.clusters} below), as references to
\autoref{tab.1}.
\endroster
The rest of the data is common to
Tables~\ref{tab.1}--\ref{tab.3.2};
they apply to a very general member $X\in\famX$ of the
respective stratum:
\roster*
\item
the numbers of lines and conics on~$X$, in the same form as in
\autoref{tab.Mukai},
\item
the order of the group $G:=\Aut\Fn(X,h)$;
if $N:=\NS(X)$ is \emph{not} generated by lines and conics, it is shown in the
form $\ls|\OG_h(N)|\cdot[G:\OG_h(N)]$,
\item
the index $i_\dm:=[G:G_\dm]$ of the setwise stabilizer
$G_\dm:=\stab\dm\subset G$,
\item
the group (as a reference to \autoref{tab.groups})~$G_\Go$ of symplectic
automorphisms of~$X${\3, as well as}
the index $[\Aut_hX:G_\Go]$, if greater
than~$1$, as a superscript,
\item
the determinant $\ls|\det\NS(X)|=\ls|\det T(X)|$ and the index
$[\NS(X):\Fano(\graph)]$, if greater than~$1$, as a superscript
(see also \autoref{rem.tables}),
\item
the numbers $(r,c)$ of, respectively, real components
and pairs of complex conjugate
components of the stratum, see \autoref{rem.tables}.
\endroster

\remark\label{rem.tables}
In Tables~\ref{tab.3} and~\ref{tab.3.2} listing the singular octics
{\2(in the sense of \emph{singular $K3$-surfaces},
\ie, those of the maximal Picard rank $\rho=20$)}, instead
of $\det T(X)$ we show the isomorphism classes of the
transcendental lattice~$T(X)$, each class in a separate row. The counts
$(r,c)$ are itemized accordingly.

Given a pair $(\graph,\dm)$ and a class $T\in\operatorname{genus}T(X)$,
the counts $(\br,\bc)$ for the relative stratum
\smash{$\famXdm_T(\graph,\dm)$}
may differ from the respective counts $(r,c)$ for $\famX_T(\graph)$.
If this is the case, the counts are shown in the form $(r,c)\NE(\br,\bc)$.
\endremark

\subsection{Restrictions on extra lines and conics}\label{s.restrictions}
We start with a few further (\ie, beyond those found in \autoref{s.supports})
restrictions on the
supports of an extra line or conic~$u$. Note that the statement and proof of
\autoref{lem.supports}, as well as those of \autoref{lem.char} concerning the
case $\supp_1u\in\CO_*$, are valid for any degree $h^2\in4\Z^+$.

\lemma[see \cite{degt:4Kummer}]\label{lem.supports}
Let $u\notin\lattice$ be an extra conic \rom(line\rom), and let
\[*
\gu:=\supp_1u,\ p:=\ls|\gu|,\quad\text{and}\quad
\ggu:=\supp_2u,\ q:=\ls|\ggu|.
\]
Then, for any pair $\gv,\ggv\subset\dm$ such that
\[*
\gv\in\cl\gu_p,\quad\ggv\subset\dm\sminus\gv,\quad\ls|\ggv|=q,
\]
there is a conic \rom(resp.\ line\rom) $v\in\lattice\(u)$ such that
$\supp_1v=\gv$ and $\supp_2v=\ggv$.
\endlemma

\proof
For completeness, we cite the proof found in~\cite{degt:4Kummer}.
A set $\gv$ as in the statement has the form $\gv=\gu\sd\go$ for some
$\go\in\CC_*$ such that $2\ls|\go\cap\gu|=\ls|\go|$.
Let $\gs_+:=\go\cap\ggu$ and
pick  $\gs_-\subset\go\cap\gu$ so that $\ls|\gs_-|=\ls|\gs_+|$.
Then, the vector
\[*
w:=u+\ds\go\gu+\gs_+-\gs_-\in\lattice\(u)
\]
has $\supp_1w=\gu\sd\go$ and
$\ggw:=\supp_2w=\ggu\sd(\gs_+\cup\gs_-)$,
so that $\ls|\ggw|=\ls|\ggu|$.
There remains to let
\[*
v:=w+(\ggv\sminus\ggw)-(\ggw\sminus\ggv).
\qedhere
\]
\endproof

%
%
%
%
%

\lemma[\cf. \cite{degt:4Kummer}]\label{lem.char}
If $u\notin\lattice$ is an extra conic and $\gu:=\supp_1u\in\KC\cup\KK$,
then any geometric extension of the lattice $\lattice\(u)$ is generated by
lines over $\lattice$.
\endlemma

\proof
Assuming the contrary, let $\ggu:=\supp_2u$ and consider the vector
\[*
\hu:=
\begin{cases}
  u-\ds\gu\varnothing+\ggu, & \mbox{if $\gu\in\KC$},\\
  \hh-u-\ds{\bar\gu}{\ggu}, & \mbox{if $\gu\in\KK$}.
\end{cases}
\]
We have $\hu\in\Sperp$, see \eqref{eq.T}, and, respectively,
\[*
\alignedat3
\hu^2&=\phantom{-}\tfrac12\ls|\gu|+2\ls|\ggu|-2,&
 \quad\hu\cdot h&=\phantom{0}2+\ls|\gu|+2\ls|\ggu|
 &\quad&\text{if $\gu\in\KC$},\\
\hu^2&=-\tfrac12\ls|\gu|+\tfrac14h^2+4,&
 \quad\hu\cdot h&=14-\ls|\gu|-2\ls|\ggu|+\tfrac12h^2
 &\quad&\text{if $\gu\in\KK$}.
\endalignedat
\]
In view of~\eqref{eq.Sperp}, the presence of this vector $\hu\in\Sperp$ rules out
the patterns
$\cn p-0$, $p=0$, $8$, $16$. The few remaining cases
(see \autoref{tab.Sylvester}) are considered
below.

\emph{The patterns $\cn12-q$, $q=2,4$}:\enspace
we have $\hu^2=0$ and $\hu\cdot h=\pm2$, \ie, $\hu$ is a $2$-isotropic vector,
see \autoref{def.admissible}\iref{i.isotropic}.

\emph{The patterns $\cn0-1$ and $\cn12-q$, $q=0,1$}:\enspace
we have $\hu^2=0$ and $\hu\cdot h=6$ or $4$. Therefore,
by \autoref{lem.kernels}, any geometric extension of $\lattice\(u)$ must
contain a vector of the form $v:=-\frac12\hu-\ds{\gs}\varnothing$
for some $\gs\in\CC_0\cup\CC_4$.
If $\gs\in\CC_0$, \ie, $\gs=\varnothing$, then $\hu$ is divisible
by~$2$; due to~\eqref{eq.Sperp}, this is only possible if $\hu\cdot h=4$,
making $\frac12\hu$ a $2$-isotropic vector, see
\autoref{def.admissible}\iref{i.isotropic}. Otherwise, if $\gs\in\CO_4$, we
obtain
\[*
v^2=-2,\quad\text{$v\cdot h=1$ or~$2$},\quad\text{$v\cdot e=-1$ for each $e\in\gs$},
\]
resulting in a missing conic, see \autoref{def.admissible}\iref{i.missing},
or exceptional divisor $v-e$, see
\autoref{def.admissible}\iref{i.exceptional}, respectively.

\emph{The pattern $\cn4-0$}:\enspace
we have $\hu^2=4$ and $\hu\cdot h=14$; by \autoref{lem.kernels}, any
geometric extension of $\lattice\(u)$ must contain a line of the form
$\frac12\hu+\ds{\gs}\varnothing$, $\gs\in\CC_6$.
Observe that, in fact, this is the only case where the lattice $\lattice\(u)$
as in the statement does admit a geometric extension, \cf.
\autoref{lem.line.6} below.
\endproof

\lemma\label{lem.supp2=0}
Let $u\notin\lattice$ be an extra conic and
assume that
$\ggu:=\supp_2u\ne\varnothing$.
Then the lattice $\lattice\(u)$ has no geometric extensions.
\endlemma

\proof
According to Tables~\ref{tab.orbits},~\ref{tab.Sylvester}
and \autoref{lem.char}, we can assume
that
\[*
\gu:=\supp_1u\in\CC_{12}\sminus\KK;
\]
hence, there is a set $\gk\in\CK_4$ such that $\ls|\gk\cap\gu|=2$.
Using \autoref{lem.supports}, we can change the set~$\ggu$ so that
$\ls|\gk\cap\ggu|\ge\min\{2,\ls|\ggu|\}$.
Pick a singleton $\gs\subset\gk$ as follows:
\roster*
\item
$\gs\subset\gk\sminus(\gu\cup\ggu)$ if $\ls|\ggu|=1$, or
\item
$\gs\subset\gk\cap\gu$ if $\ls|\ggu|\ge2$.
\endroster
Then, for the \BB-conic $v:=\hh+\ds\gk\gs$, we have $v\cdot u=-1$ and, hence,
$u-v$ is an exceptional divisor, see
\autoref{def.admissible}\iref{i.exceptional}.
\endproof

\lemma\label{lem.class.2}
Let $u\notin\lattice$ be an extra conic,
$\gu:=\supp_1u$, and
$p:=\ls|\gu|$.
Then, for any {\3given} set $\gw\in\Cl\gu_{16-p}\sminus\cl\gu$,
there is a conic $w\in\lattice\(u)$ such
that $\supp_1w=\gw$.
\endlemma

\proof
Any set $\gw$ as in the statement is of the form $\overline{\gv\sd\gs}$,
where $\gv:=\supp_1v\in\cl\gu_p$ for an appropriate vector~$v$ given by
\autoref{lem.supports} and $\gs\in\CK_4$, $\ls|\gs\cap\gv|=2$.
Besides, by \autoref{lem.supp2=0} we can assume that
$\supp_2u=\supp_2v=\varnothing$.
Then, it is immediate that the conic
$w:=\hh-\ds\gs\gv-v$ is as required.
\endproof

\lemma\label{lem.line.8}
Let $u\notin\lattice$ be an extra line and
$\gu:=\supp_1u\in\CC_8$.
Then, the lattice $\lattice\(u)$ is not admissible.
\endlemma

\proof
There exists a subset $\gr\in\CK_{12}$ such that $\ls|\gr\cap\gu|=7$; then, it is
immediate that
$-\hh+\ds\gr\gu+2u$ is an exceptional divisor, see
\autoref{def.admissible}\iref{i.exceptional}.
\endproof

\lemma\label{lem.line.6}
Let $u\notin\lattice$ be an extra line and
$\gu:=\supp_1u\in\CC_6$.
Then\rom:
\roster
\item\label{i.BB}
all sixteen \BB-conics, see \eqref{eq.BB}, are reducible in $\lattice\(u)$\rom;
\item\label{i.4-0}
for each $\gv\in\CK_4$, there is an irreducible conic $v\in\lattice\(u)$ with
$\supp_1v=\gv$\rom;
\item\label{i.12-0}
for each $\gv\in\Cl\gu_{12}$, there is a line $v\in\lattice\(u)$ with
$\supp_1v=\gv$.
\endroster
\endlemma

\proof
For statement~\iref{i.BB}, observe that, for each pair
$\gs\subset\gk\in\CK_4$ as in~\eqref{eq.BB}, there is
$\gw\in\cl\gu_6$ such that $\gw\cap\gk=\gk\sminus\gs$; then,
$w\cdot k=-1$,
where $w\in\lattice\(u)$ is the line with $\supp_1w=\gw$ given by
\autoref{lem.supports} and
$k=\hh+\ds\gk\gs$ is the \BB-conic~\eqref{eq.BB}.

For each pair $\gk$, $\gw$ as above, the line $v:=k-w$ has support
$\gv:=\overline{\gw\sd\gk}\in\Cl\gu_{12}$, and all lines as in
statement~\iref{i.12-0} can be obtained in this way.

Finally, the four extra conics as in statement~\iref{i.4-0} are
\[*
\hh-\ds\gr\gw-(\gv\sminus\gr)-2w,
\]
where $\gw$ and $w$ are as above and
$\gr\in\CK_{12}$, $\ls|\gr\cap\gw|=3$; \cf. the last case $\cn4-0$ in the
proof of \autoref{lem.char}.
\endproof

\subsection{Proof of \autoref{th.codim.1}}\label{proof.codim.1}
According to Tables~\ref{tab.orbits}, \ref{tab.Sylvester} and
Lemmas~\ref{lem.char}, \ref{lem.supp2=0}, \ref{lem.line.8}, there are but five (pairs of)
patterns that need to be considered:
\[*
\cn4-0,\cn12-0;\ \cn6-0,\cn10-0;\ \cn8-0\quad\text{or}\quad
\ln4-0;\ \ln6-0,\ln10-0.
\]
Here, two patterns constitute a pair, \eg, $\cn4-0,\cn12-0$, if they result
in identical $1$-vector extensions: in the example, the extension
$\lattice\(u)$ by a vector with pattern $\cn4-0$ contains one with pattern
$\cn12-0$ (see \autoref{lem.class.2} or, for lines, \autoref{lem.line.6})
and \latin{vice versa}.

Furthermore, \autoref{lem.char} asserts that
$\gu:=\supp_1u\notin(\CC_*\cup\CK_*)$:
the case $\gu\in\CK_4$ can be ignored
as the lattice $\lattice\(u)$ itself is not geometric whereas
any geometric extension thereof is generated by lines, \viz. the pair of
patterns $\ln6-0,\ln10-0$. Obviously, the $\G$-isomorphism class of
$\lattice\(u)$ depends only on the $\G$-orbit of $\gu$; by
\autoref{lem.supports}, this can further be replaced with the $\G$-orbit of
$\cl{\gu}$. Hence, referring to
\autoref{tab.orbits} and parity condition~\eqref{eq.parity},
we conclude that each of
the five (pairs of) patterns above results in a single $\G$-isomorphism
class of extensions.
{\2Now, a straightforward
computation based on \autoref{s.computation}
shows that}
\roster*
\item
each of the five
lattices $N:=\lattice\(u)$ obtained in this way is geometric,
\item
there are no
proper geometric finite index extensions $N'\supset N$, and
\item
each lattice $N\supset\lattice\supset\dm$ admits a unique
$\OG_h(N,\dm)$-isomorphism class of primitive isometries
$N\into\L$ (see \autoref{s.components}).
\endroster
Thus, there are five strata, each consisting of a single real
component (see \autoref{s.components}),
and using~\eqref{eq.Fn} one can
{\2compute the Fano graphs and, in particular,} show that,
in addition to~$\dm$ and \BB-conics~\eqref{eq.BB},
the lines and conics in~$N$ are
{\2exactly}
those given by
Lemmas~\ref{lem.supports}, \ref{lem.class.2}, and~\ref{lem.line.6}.
The precise counts are given in \autoref{tab.1}.
\qed

\subsection{Clusters}\label{s.clusters}
The discriminant $\discr\lattice$ has $2$- and $5$-torsion:
\[*
\discr_2\lattice\cong
 \bmatrix0&\frac12\\\frac12&0\endbmatrix\oplus\bmatrix0&\frac12\\\frac12&0\endbmatrix
  \oplus\bigl[\tfrac58\bigr],\qquad
\discr_5\lattice\cong\bigl[\tfrac85\bigr].
\]
The groups $2\discr_2\lattice\cong\Z/4$ and $\discr_5\lattice\cong\Z/5$
{\2(where $\discr_p:=\Z_p\otimes\discr$)}
have distinguished generators
$\eta_2:=\frac14\th$ and $\eta_5:=\frac15\th$,
respectively, see~\eqref{eq.th}.

Consider a geometric extension $N\supset\lattice$.
Following~\cite{degt:4Kummer}, define a \emph{cluster} in~$N$ as a collection
of all lines and conics $u\in N$  sent to the same point
of the projective space $\Bbb{P}((N/\lattice)\otimes\Q)$. Consider also the
canonical homomorphism
\[*
\Gd=\Gd_2\oplus\Gd_5\:N\to\lattice\dual\to\discr\lattice=\lattice\dual/\lattice.
\]
Directly by the definition, the image $\Gd(C)$ of each cluster $C\subset N$
generates a cyclic subgroup in $\discr\lattice$. More precisely, since each
cluster is contained in a $1$-vector extension, \autoref{th.codim.1} and
Lemmas~\ref{lem.supports}, \ref{lem.class.2}, \ref{lem.line.6}
used in its proof imply that the image of each cluster consists of
\roster*
\item
a single element $\Ga$, as in stratum~\strat1 in \autoref{tab.1},
or
\item
a pair of elements $\pm\Ga$, as in strata~\strat3, \strat4, \strat5, or
\item
a pair $\pm\Ga$ and common element $2\Ga=\eta_2\oplus2\eta_5$,
as in stratum~\strat2.
\endroster
The generating images $\Gd(u)=\Gd_2(u)\oplus\Gd_5(u)$ are shown in
\autoref{tab.1},
{\3as}
the square $\Gd_2^2=r/s\bmod2\Z$ (where $s$ is the order
of~$\Gd_2$) and coefficient of $\Gd_5$ in the basis~$\eta_5$. Computing the
orbits of the $\G$-action on $\discr\lattice$, we conclude that,
with the extra restriction that
\[*
\aligned
&\text{$\Gd_2(u)\cdot\eta_2=\tfrac14(\epsilon+p)\bmod\Z$ for $u$ with pattern
 $\cn p-0$ ($\Ge=2$) or $\ln p-0$ ($\Ge=1)$},\\
&\text{$\Gd_2(u)\ne\pm\eta_2$ unless $u$ is a non-generating conic of pattern
$\cn4-0$ in stratum \strat2},
\endaligned
\]
these data determine the $\G$-orbit of $\Gd(u)$. On the other hand, by
comparison to \autoref{tab.orbits}, the vector $\Gd(u)$ determines
$\cl{\supp_1u}$ and, hence,
the extension
$\lattice\(u)$.

\section{Strata of higher codimension}\label{S.codim.large}

In this section, we complete the proofs of
the principal results of the paper by analyzing the
double and triple (self-)intersections of the five strata found in
\autoref{S.codim.1}.

\theorem\label{th.codim.2}
The space $\famB$ has $15$ irreducible equiconical strata of codimension~$2$,
see \autoref{tab.2}.
Each stratum consists of a single real component\rom; one of the absolute
strata
splits into two relative ones \rom(prefixed with \rsame21 in
\autoref{tab.2}\rom).
\table
\caption{Strata of codimension~$2$ \rm(see \autoref{th.codim.2})}\label{tab.2}
\centerline{\vbox{\small\rm\offinterlineskip
\tabdefs
\def\SAME{2}%
\tabskip0pt plus 1fil
\halign{\strut\quad$#$\hss\sep\tabskip0pt&
 \strut\sep\hss$#$\hss\sep&
 \strut\sep\hss$#$\hss\sep&
 \strut\sep\hss$#$\hss\sep&
 \strut\sep\hss$#$\hss\sep&
 \strut\sep\hss$#$\hss\sep&
 \strut\sep\hss$#$\hss\sep&
 \strut\sep$#$\hss\sep\tabskip0pt plus 1fil\cr
\noalign{\hrule\kern3pt}
\text{Clusters}&
 \text{Lines}&\text{Conics}&\ls|G|&i_\dm&
 G_\Go&\ls|{\det}|&\omit\hss$(r,c)$\hss\cr
\noalign{\kern1pt\hrule\kern3pt}
\s1,\s1&8&32&384&2&\(16,14)\^1&240&(1,0)\EQ(1,0)\cr
\s1,\s1,\s5&8&8\+72&256&2&\(16,14)\^2&160&(1,0)\EQ(1,0)\cr
\s1,\s2,\s4&24&32\+36&192&2&\(16,14)\^1&80&(1,0)\EQ(1,0)\cr
\noalign{\kern1pt\hrule\kern3pt}
\same{1} \s1,\s3&4&40&64&1&\(16,14)\^1&320&(1,0)\EQ(1,0)\cr
\same{1} \s1,\s3&4&40&64&1&\(16,14)\^1&320&(1,0)\EQ(1,0)\cr
\s1,\s4&4&64&384&4&\(16,14)\^1&240&(1,0)\EQ(1,0)\cr
\s2,\s3&20&16\+28&64&1&\(16,14)\^1&128&(1,0)\EQ(1,0)\cr
\noalign{\kern1pt\hrule\kern3pt}
\s3,\s3&&48&256&2&\(16,14)\^1&416&(1,0)\EQ(1,0)\cr
\s3,\s3&&48&512\cdot16&2&\(16,14)\^1&512\^2&(1,0)\EQ(1,0)\cr
\s3,\s3&&48&512&2&\(16,14)\^1&512&(1,0)\EQ(1,0)\cr
\s3,\s3,\s4&&80&512&4&\(16,14)\^1&288&(1,0)\EQ(1,0)\cr
\s3,\s4&&72&512&4&\(16,14)\^1&320&(1,0)\EQ(1,0)\cr
\s3,\s5&&88&256&2&\(16,14)\^2&288&(1,0)\EQ(1,0)\cr
\s4,\s4&&96&2304&6&\(48,50)\^1&224&(1,0)\EQ(1,0)\cr
\s4,\s5&&112&1024&4&\(16,14)\^2&192&(1,0)\EQ(1,0)\cr
\maxr1\s5,\s5&&128&1024&2&\(32,27)\^2&160&(1,0)\EQ(1,0)\cr
\noalign{\kern1pt\hrule}
\crcr}}}
\endtable
\endtheorem

In a stratum of codimension~$3$, each octic~$X$ is a
so-called \emph{singular $K3$-surface}
($\rank\NS(X)=20$ is maximal); hence, $X$ is \emph{rigid}, \ie,
$X$ is projectively equivalent to any equiconical deformation thereof. In
other words, modulo the group $\PGL(\C,6)$,
the union of the codimension~$3$
strata is a finite collection of points,
and it is these points that are
listed in Tables~\ref{tab.3} and~\ref{tab.3.2}. {\2(In particular, this list
also proves the finiteness of the moduli space; we refrain from discussing the
general algebra-geometric philosophy behind this phenomenon.)}

\theorem\label{th.codim.3}
All equiconically rigid Barth--Bauer octics are listed in
Tables~\ref{tab.3},~\ref{tab.3.2}\rom;
\table
\caption{Rigid octics with $>80$ conics \rm(see \autoref{th.codim.3})}\label{tab.3}
\centerline{\vbox{\small\rm\offinterlineskip
\tabdefs
\def\SAME{3}%
\tabskip0pt plus 1fil
\halign{\strut\quad$#$\hss\sep\tabskip0pt&
 \strut\sep\hss$#$\hss\sep&
 \strut\sep\hss$#$\hss\sep&
 \strut\sep\hss$#$\hss\sep&
 \strut\sep\hss$#$\hss\sep&
 \strut\sep\hss$#$\hss\sep&
 \strut\sep\hss$#$\hss\sep&
 \strut\sep$#$\hss\sep\tabskip0pt plus 1fil\cr
\noalign{\hrule\kern3pt}
\text{Clusters}&\text{Lines}&\text{Conics}&\ls|G|&i_\dm&
 G_\Go&T&\omit\hss$(r,c)$\hss\cr
\noalign{\kern1pt\hrule\kern3pt}
\s5,\s5,\s5&&176&15360&10&\(960,11357)\^2&[8,4,12]&(1,0)\EQ(1,0)\cr
\s4,\s5,\s5&&160&3072&12&\(192,1493)\^2&[4,0,24]&(1,0)\EQ(1,0)\cr
\s3,\s5,\s5&&136&512&2&\(32,27)\^2&[4,0,36]&(1,0)\EQ(1,0)\cr
\s1,\s1,\s1,\s1,\s5,\s5&16&32\+96&256&2&\(32,27)\^2&[4,2,16]&(1,0)\EQ(1,0)\cr
\s3,\s3,\s3,\s4,\s4&&120&384&6&\(48,50)\^1&[8,4,20]&(1,0)\NE(0,1)\cr
\s3,\s4,\s5&&120&256&4&\(16,14)\^2&[8,0,20]&(1,0)\NE(2,0)\cr
\s1,\s1,\s4,\s5&8&8\+104&256&4&\(16,14)\^2&[4,0,24]&(1,0)\NE(0,1)\cr
\maxr2\s1,\s1,\s2,\s4,\s4&28&48\+52&288&3&\(48,50)\^1&[4,2,12]&(1,0)\EQ(1,0)\cr
\s1,\s4,\s4&4&96&576&6&\(48,50)\^1&[4,2,36]&(1,0)\NE(2,0)\cr
\s3,\s3,\s5&&96&256&2&\(16,14)\^2&[8,4,28]&(1,0)\EQ(1,0)\cr
\s3,\s3,\s5&&96&256&2&\(16,14)\^2&[8,0,32]&(1,0)\EQ(1,0)\cr
\s3,\s3,\s5&&96&256&2&\(16,14)\^2&[8,0,32]&(1,0)\EQ(1,0)\cr
\s3,\s3,\s3,\s3,\s4&&96&256&4&\(32,27)\^1&[8,0,24]&(1,0)\NE(0,1)\cr
\s3,\s3,\s3,\s4&&88&128&4&\(16,14)\^1&[8,4,32]&(1,0)\NE(0,2)\cr
\s1,\s1,\s3,\s5&8&8\+80&32&2&\(16,14)\^2&[4,2,32]&(0,1)\EQ(0,1)\cr
&&&&&&[8,2,16]&(0,2)\EQ(0,2)\cr
\s1,\s2,\s3,\s3,\s4&24&32\+52&64&2&\(16,14)\^1&[8,2,8]&(1,0)\NE(2,0)\cr
\noalign{\kern1pt\hrule}
\crcr}}}
\endtable
\table
\caption{Other rigid octics \rm(see \autoref{th.codim.3})}\label{tab.3.2}
\centerline{\vbox{\small\rm\offinterlineskip
\tabdefs
\def\SAME{3}%
\tabskip0pt plus 1fil
\halign{\strut\quad$#$\hss\sep\tabskip0pt&
 \strut\sep\hss$#$\hss\sep&
 \strut\sep\hss$#$\hss\sep&
 \strut\sep\hss$#$\hss\sep&
 \strut\sep\hss$#$\hss\sep&
 \strut\sep\hss$#$\hss\sep&
 \strut\sep\hss$#$\hss\sep&
 \strut\sep$#$\hss\sep\tabskip0pt plus 1fil\cr
\noalign{\hrule\kern3pt}
\text{Clusters}&\text{Lines}&\text{Conics}&\ls|G|&i_\dm&
 G_\Go&T&\omit\hss$(r,c)$\hss\cr
\noalign{\kern1pt\hrule\kern3pt}
\s1,\s3,\s3,\s4&4&80&64&4&\(16,14)\^1&[8,2,20]&(0,1)\NE(0,4)\cr
\s3,\s3,\s4&&80&256&4&\(32,27)\^1&[12,4,20]&(0,1)\NE(0,2)\cr
\s1,\s1,\s1,\s1,\s5&16&16\+64&512&2&\(32,27)\^2&[8,4,12]&(1,0)\EQ(1,0)\cr
\s1,\s2,\s3,\s4&24&32\+44&64&2&\(16,14)\^1&[4,0,16]&(1,0)\NE(0,1)\cr
\same{1} \s1,\s3,\s4&4&72&64&2&\(16,14)\^1&[4,0,44]&(1,0)\NE(0,1)\cr
&&&&&&[12,4,16]&(0,1)\NE(0,2)\cr
\same{1} \s1,\s3,\s4&4&72&64&2&\(16,14)\^1&[4,0,44]&(1,0)\NE(0,1)\cr
&&&&&&[12,4,16]&(0,1)\NE(0,2)\cr
\s1,\s1,\s4&8&64&256&4&\(32,27)\^1&[12,0,12]&(0,1)\NE(0,2)\cr
\s3,\s3,\s3,\s3&&64&256&2&\(32,27)\^1&[8,0,32]&(1,0)\NE(2,0)\cr
\s3,\s3,\s3&&56&384&2&\(48,50)\^1&[4,0,68]&(1,0)\NE(2,0)\cr
&&&&&&[8,4,36]&(1,0)\NE(0,1)\cr
\s3,\s3,\s3&&56&64&2&\(16,14)\^1&[8,4,48]&(1,0)\NE(0,1)\cr
&&&&&&[16,4,24]&(0,1)\NE(0,2)\cr
\s2,\s3,\s3&20&16\+36&64&1&\(16,14)\^1&[8,4,16]&(1,0)\EQ(1,0)\cr
\s2,\s3,\s3&20&16\+36&64&1&\(16,14)\^1&[8,4,16]&(1,0)\EQ(1,0)\cr
\s2,\s3,\s3&20&16\+36&32&1&\(16,14)\^1&[4,2,24]&(1,0)\EQ(1,0)\cr
&&&&&&[8,2,12]&(0,1)\EQ(0,1)\cr
\s1,\s1,\s3,\s3&8&48&64&2&\(16,14)\^1&[8,2,20]&(0,1)\NE(0,2)\cr
\s1,\s3,\s3&4&48&64&2&\(16,14)\^1&[16,0,16]&(0,1)\NE(0,2)\cr
\s1,\s3,\s3&4&48&64&2&\(16,14)\^1&[16,0,16]&(0,1)\NE(0,2)\cr
\same{2} \s1,\s3,\s3&4&48&64&1&\(16,14)\^1&[8,4,32]&(2,0)\EQ(2,0)\cr
\same{2} \s1,\s3,\s3&4&48&64&1&\(16,14)\^1&[8,4,32]&(2,0)\EQ(2,0)\cr
\same{3} \s1,\s3,\s3&4&48&64&1&\(16,14)\^1&[8,4,32]&(2,0)\EQ(2,0)\cr
\same{3} \s1,\s3,\s3&4&48&64&1&\(16,14)\^1&[8,4,32]&(2,0)\EQ(2,0)\cr
\same{4} \s1,\s3,\s3&4&48&32&1&\(16,14)\^1&[4,2,56]&(2,0)\EQ(2,0)\cr
&&&&&&[16,6,16]&(0,1)\EQ(0,1)\cr
\same{4} \s1,\s3,\s3&4&48&32&1&\(16,14)\^1&[4,2,56]&(2,0)\EQ(2,0)\cr
&&&&&&[16,6,16]&(0,1)\EQ(0,1)\cr
\same{5} \s1,\s1,\s3&8&40&64&1&\(16,14)\^1&[4,0,44]&(1,0)\EQ(1,0)\cr
&&&&&&[12,4,16]&(0,1)\EQ(0,1)\cr
\same{5} \s1,\s1,\s3&8&40&64&1&\(16,14)\^1&[4,0,44]&(1,0)\EQ(1,0)\cr
&&&&&&[12,4,16]&(0,1)\EQ(0,1)\cr
\s1,\s1,\s1&12&32&576&2&\(48,50)\^1&[4,2,36]&(1,0)\NE(2,0)\cr
\noalign{\kern1pt\hrule}
\crcr}}}
\endtable
altogether, there are
\roster*
\item
$36$ isomorphism classes of abstract Fano graphs~$\graph$,
\item
$41$ isomorphism classes of pairs $(\graph,\dm)$,
\item
$33$ real and $14$ pairs of complex conjugate octics~$X$, and
\item
$38$ real and $38$ pairs of complex conjugate pairs~$(X,\dm)$.
\endroster
\endtheorem

\subsection{Proof of Theorems~\ref{th.codim.2} and~\ref{th.codim.3}}\label{proof.codim.2.3}
We use the approach of~\cite[\S3]{degt:4Kummer}.

For \autoref{th.codim.2}, we consider all corank~$2$ extensions $\lattice\(u,v)$
by a pair of vectors, each as in \autoref{tab.1}; an extra piece of
data is the product $u\cdot v$, which must satisfy \autoref{lem.intersection}.
(We adopt Convention~3.9 in~\cite{degt:4Kummer} and assume that the
generating set has the maximal number of lines; then, we can also assume
that all generating conics are irreducible and, hence, $u\cdot v\ge0$.)
The vast majority of possibilities are ruled out by the Hodge index theorem,
as in \autoref{s.supports},
leaving but $30$ $\G$-orbits of triples
$(\cl\gu,\cl\gv,u\cdot v)$.
Each triple is analyzed in the spirit of \autoref{S.codim.1}, and
only $20$ of them admit a geometric finite index
extension (which is always trivial).
There remains to observe that some of the lattices obtained are
isomorphic: in fact, each geometric lattice $\lattice\(u,v)$
is generated over~$\lattice$ by appropriate representatives of any pair of
clusters contained in $\lattice\(u,v)$.

\autoref{th.codim.3} is proved similarly, by extending one of the
$16$ geometric lattices
$\lattice\(u,v)$ given by \autoref{th.codim.2} by a third extra line or
conic~$w$.
\qed

\subsection{Proof of \autoref{th.main}}\label{proof.main}
The bound $\ls|\Fn_2X|\le176$ and the uniqueness of the Barth--Bauer octic
$X_{176}$ at which this bound is attained are given by Theorems~\ref{th.codim.1},
\ref{th.codim.2}, \ref{th.codim.3}. Furthermore,
$X_{176}$ admits a faithful projective symplectic action of the Mukai
group~$M_{20}$ (see~\cite{Mukai};
\#\symplectic(960,11357) in \autoref{tab.groups}).
On the other hand, according to \cite[Corollary~7.3]{degt:4Kummer}
(see also~\cite{Bonnafe.Sarti},
where a slightly stronger assumption is used),
this property characterizes a unique
octic $K3$-surface $X\subset\Cp5$.
The
defining equations cited in \autoref{th.main} are found
in~\cite{Bonnafe.Sarti}.
\qed

\example\label{ex.160}
It is remarkable that the only Barth--Bauer octic~$X_{160}$
realizing the next largest
number $160$ of conics (the second row in \autoref{tab.3})
is also characterized by the presence of a faithful projective symplectic
action of a Mukai group, this time $T_{192}$ (\#\symplectic(192,1493) in
\autoref{tab.groups}). The uniqueness of a $T_{192}$-octic in $\Cp5$
is easily proved
similar to \cite[\S\,7.1]{degt:4Kummer}.

First, the N\'{e}ron--Severi lattice~$S$ of a very general (non-algebraic)
$K3$-surface with a faithful symplectic $T_{192}$-action (\cf.
\cite{Hashimoto}) can be found as $h^\perp\subset\NS(X_{160})$. One has
\[*
\discr_2S=\bigl[\tfrac54\bigr]\oplus
    \bigl[\tfrac54\bigr]\oplus\bigl[\tfrac54\bigr],\qquad
\discr_3S=\bigl[\tfrac43\bigr],
\]
and the image of the natural homomorphism
$\Aut(\Fn X_{160})\into\OG(S)\to\Aut(\discr S)$
is an index~$12$ subgroup preserving one of the $12$ vectors~$\Ga_i$ of square
$\tfrac32\bmod2\Z$.

On the other hand, each of the twelve vectors~$\Ga_i$ as above gives rise to
an index~$4$ extension of $S\oplus\Z h$, which is the N\'{e}ron--Severi lattice
of a Barth--Bauer octic with $160$ conics. By \autoref{th.codim.3}, we
conclude that all these extensions are isomorphic;
hence, all $12$ vectors
constitute a single $\OG(S)$-orbit and the natural homomorphism
$\OG(S)\to\Aut(\discr S)$ is surjective.

From the last statement, using the techniques of~\cite{Nikulin:forms}
{\2outlined in \autoref{s.computation}\iref{i.T.genus},~\iref{i.T.cosets}}
and the uniqueness of
\[*
S^\perp\cong\bmatrix
  4&  0&  0\\
  0&  8&  4\\
  0&  4&  8
\endbmatrix
\]
in its genus, we conclude that there is a single $\OG(S)$-equivalence class
of primitive isometries $S\into\L$; furthermore,
{\2in view of~\eqref{eq.gh.ext},}
any element of
$\OG(S^\perp)$ extends to an autoisometry of~$\L$.
Since the group $\OGplus(S^\perp)$ acts transitively on the six square~$8$ vectors in
$S^\perp$, the uniqueness of a $T_{192}$-octic surface follows, \cf.
\autoref{s.components}.
\endexample

\remark\label{rem.X64}
The same argument shows that there is a unique $T_{192}$-quartic in~$\Cp3$.
It is the famous Schur~\cite{Schur:quartics} quartic~$X_{64}$ maximizing the
number of lines: it has $64$ lines and
$\text{$576$ reducible}+\text{$144$ irreducible}=720$ conics.
\endremark

\subsection{Proof of \autoref{th.lines}}\label{proof.lines}
The bound on the number of lines is explicitly stated in~\cite{degt:lines}.
To estimate the number of reducible conics (\ie, pairs of intersecting
lines), recall the bound
\[*
\operatorname{val}v\le\begin{cases}
           7, & \mbox{if $X$ is a triquadric}, \\
           8, & \mbox{if $X$ is a special octic}
         \end{cases}
\]
on the valency of a line in the graph $\Fn_1X$, see
\cite[Proposition~2.12]{degt:lines}. It follows that the number of reducible
conics does not exceed
\[*
\begin{cases}
  30\cdot7/2=105, & \mbox{if $X$ is a triquadric and $\ls|\Fn_1X|\le30$}, \\
  26\cdot8/2=104, & \mbox{if $X$ is special and $\ls|\Fn_1X|\le26$}.
\end{cases}
\]
On the other hand, the Fano graphs of the triquadrics with more than $30$
lines and special octics with more than $26$ lines are listed
in~\cite{degt:lines} (see Theorems~1.2 and~1.4 respectively), and the number
of reducible conics in these graphs is easily computed: the maximum is $112$,
attained at a unique triquadric{\3, \viz. the one denoted by $\Theta'_{36}$
in~\cite{degt:lines}.}
\qed

\subsection{Proof of \autoref{th.real}}\label{proof.real}
As explained in~\cite{degt:4Kummer},
{\em an equiconical stratum of Barth--Bauer
octics contains a real octic with all lines and conics real if and only if
the respective generic transcendental lattice has a direct
summand isomorphic to~$\bU(2)$.} In particular, this stratum must have
codimension at most~$2$. On the other hand, according to
Theorems~\ref{th.codim.1} and~\ref{th.codim.2}, the maximal number of conics
on a Barth--Bauer octic of Picard rank $\rho\le19$ is $128$ (see the line
marked with a \maxreal1 in \autoref{tab.2}), the typical transcendental
lattice being $T\cong\bU(2)\oplus[40]$, as required.

To show that this is the maximum, we have to consider singular octics
given by \autoref{th.codim.3} and Tables~\ref{tab.3},~\ref{tab.3.2} and, for
each such octic~$X$, compute the actions~$c_*$ induced on $\NS(X)$ by
all possible real structures $c\:X\to X$.
{\2Arithmetically, we consider involutive elements
$c_\graph\in\Aut\graph$, $\graph:=\Fn X$, with the following properties:
\roster
\item
$c_\graph$ extends to $\NS(X)$, see~\eqref{eq.g.ext}:
this requirement is redundant as we have $\NS(X)=\Fano(\graph)$ in all cases;
\item
there is an involution $c_T\in\OG(T)\sminus\OG^+(T)$ such that
$c_\graph\oplus c_T$ extends to~$\L$, see~\eqref{eq.gh.ext};
\endroster
then, $-(c_\graph\oplus c_T)$ is induced by a real structure.
(Recall that a real structure reverses the orientation of algebraic curves
and takes $H^{2,0}$ to $H^{0,2}$, see \cite{degt:4Kummer} for
details.)}
This {\2\GAP~\cite{GAP4} aided} computation gives us at most $56$ real conics.
{\2In fact, all maximal configurations correspond}
to certain real structures on the octic
$X_{176}$ introduced in \autoref{th.main}
(the first row in \autoref{tab.3}).
\qed

{
\let\.\DOTaccent
\def\cprime{$'$}
\bibliographystyle{scm}
\bibliography{degt}
}

\end{document}